\documentclass[preprint,11pt]{imsart}

\usepackage{a4wide,amsmath,amsfonts,dsfont,latexsym,amssymb,amscd, graphicx,verbatim,color}

\startlocaldefs

\newcommand{\R}{{\mathbb R}}

\newcommand{\eps}{{\varepsilon }}

\newtheorem{thm}{Theorem}[section]

\numberwithin{equation}{section}

\endlocaldefs

\begin{document}

\begin{frontmatter}

\title{Local maxima of two dependent Brownian Motions never coincide}
\runtitle{Local maxima of two dependent Brownian Motions never coincide}


\author{\fnms{Eric} \snm{Cator}\ead[label=e1]{e.a.cator@tudelft.nl}}
\address{Delft University of Technology\\
DIAM\\
Mekelweg 4\\
2628 CD Delft\\
The Netherlands\\
\printead{e1}}
\affiliation{Delft University of Technology}

\runauthor{Eric Cator}

\begin{abstract}
We consider two dependent Brownian motions with (possibly) different drift, and apply a result by le Gall on cone points of two dimensional Brownian motion to show that with probability one, there will not be a time that is a local maximum for both processes.
\end{abstract}



\end{frontmatter}

\section{Introduction}

In this short note we will consider the following problem: suppose $B_1(t)$ and $B_2(t)$ are independent two-sided standard Brownian motions. Define the processes
\[ X_1(t) = \phi_1(t) + \sigma B_1(t)\ \ {\rm and}\ \ X_2(t) = \phi_2(t) + \rho_1B_1(t) + \rho_2B_2(t).\]
Here, the drift functions $\phi_1$ and $\phi_2$ are assumed to have an $L^2$ derivative. Furthermore, the constants $\sigma$ and $\rho_2$ are assumed to be non-zero. The processes $X_1$ and $X_2$ are dependent Brownian motions with drift (in fact independent if $\rho_1=0$). We will prove the following statement:
\begin{thm}\label{thm}
Using the notations introduced above, it holds that with probability one, there does not exist $t\in \R$ such that $t$ is a local maximum for $X_1$ and for $X_2$.
\end{thm}
An application of this result can be found in \cite{Lopu} by H.P. Lopuha\"a and C. Durot. There the limiting distribution is calculated for a multiple monotone regression testing problem
\[
H_0:f_1=f_2=\cdots=f_J
\quad\text{ against }\quad
H_1: f_i\ne f_j\text{ for some }i\ne j
\]
where all $f_j$'s are decreasing. Think of the $f_j$'s as densities, regression functions or failure rates. Consider the test statistic based on comparing the isotonic estimators $\hat{f}_j$ to the pooled isotonic estimator $\hat{f}_0$ (if $H_0$ is true, all data are generated by the same $f$). Clearly, $\hat{f}_0$ is dependent of each $\hat{f}_j$. When calculating the asymptotic distribution of this test statistic, an important role is played by random variables $V_j$, which are locations of maxima of independent Brownian motions $W_j$ minus a parabola for $1\leq j\leq J$. However, the corresponding ``pooled'' variable $V_0$ is a similar location of the maximum for a Brownian motion $W_0$ minus a parabola, where $W_0$ is a weighted average of the $W_j$'s. In their analysis, Lopuha\"a and Durot need to prove that for $\eps\to 0$,
\[ {\mathbb P}(|V_0-V_j|\leq \eps) = o(1).\]
This follows directly from Theorem \ref{thm} (see page 28-29 in \cite{Lopu}).

It seems natural to try and prove Theorem \ref{thm} using path properties of one dimensional Brownian motion near a local maximum, of which many are known in the literature. However, it turned out that the most elegant way to prove Theorem \ref{thm} is to relate a coinciding local maximum to a path property of two dimensional standard Brownian motion, and apply a result by le Gall.

\section{Proof of main result}

We first restrict our time parameter $t$ to the open interval $(-T,T)$. Clearly, if we can prove for all $T>0$ that no simultaneous local maximum can exist in the interval $(-T,T)$, then the theorem follows. On the interval $(-T,T)$, define $\tilde{B}_1(t) = X_1(t)/\sigma$ and
\[ \tilde{B}_2(t) = \frac{\phi_2(t)}{\rho_2} - \frac{\rho_1\phi_1(t)}{\sigma\rho_2} + B_2(t).\]
Using the Cameron-Martin theorem it is clear that on the time interval $(-T,T)$, the law of $(\tilde{B}_1,\tilde{B}_2)$ is absolutely continuous with respect to the law of $(B_1,B_2)$. Also,
\[ (X_1(t),X_2(t)) = (\sigma \tilde{B}_1(t), \rho_1\tilde{B}_1(t) + \rho_2\tilde{B}_2(t)).\]
This means that the theorem follows if we can prove that with probability one, the processes $\sigma B_1(t)$ and $\rho_1B_1(t) + \rho_2B_2(t)$ do not have a simultaneous local maximum. Suppose $s\in \R$ is such a simultaneous maximum. Then there exists $\eta>0$ such that for all $t\in (s-\eta,s+\eta)$ we would have
\[ \sigma B_1(t)\leq \sigma B_1(s)\ \ \ \mbox{and}\ \ \ \rho_1B_1(t) + \rho_2B_2(t)\leq \rho_1B_1(s) + \rho_2B_2(s).\]
Define $p=(\sigma B_1(s), \rho_1B_1(s) + \rho_2B_2(s))$. Define $\cal C$ as the intersection of the two half-spaces:
\[ {\cal C} = \{x\in\R^2\ :\ \sigma x_1\leq p_1\}\cap\{x\in\R^2\ :\ \rho_1x_1+\rho_2x_2\leq p_2\}.\]
Then $\cal C$ is a cone with vertex $p$ and top angle $\alpha$, depending only on $\sigma, \rho_1$ and $\rho_2$, with $\alpha<\pi$, since $\sigma$ and $\rho_2$ are non-zero. Furthermore, the two dimensional Brownian motion $(B_1(t),B_2(t))$ lies inside $\cal C$ for the time interval $(s-\eta,s+\eta)$, and touches $p$ at time $s$. This makes $p$ a two-sided cone point with angle $\alpha<\pi$ for the two dimensional Brownian motion, in the sense of \cite{leGall}. However, in \cite{leGall} p.136, it is proven that with probability one, there do not exist any two-sided cone points with angle $\alpha<\pi$.\hfill $\Box$


\begin{thebibliography}{10}
\bibitem{leGall}
{\textsc{le Gall, Jean-Francois.}~(1992)}
Some properties of planar Brownian motion.
\textit{Ecole d'\'et\'e de probabbilit\'es de Saint-Flour XX - 1990.}
\textit{Lecture Notes in Mathematics}
\textbf{1527/1992}, 111--229.

\bibitem{Lopu}
{\textsc{Durot, C. and Lopuha\"a, H.P.}~(2012)}
Testing equality of functions under monotonicity constraints.
\textit{Submitted}

\end{thebibliography}
\end{document}